\documentclass[10pt]{amsart}

\usepackage[colorlinks]{hyperref}
\usepackage{color,graphicx,shortvrb}
\usepackage[latin 1]{inputenc}

\usepackage[active]{srcltx} 

\usepackage{enumerate}
\usepackage{amssymb}

\newtheorem{theorem}{Theorem}[section]

\newtheorem{lemma}[theorem]{Lemma}

\newtheorem{proposition}[theorem]{Proposition}

\def\11{\textbf{$1$}}
\def\CC{{\mathbb{C}}}

\DeclareGraphicsExtensions{.jpg,.pdf,.png,.eps}

\begin{document}

\title[On the extension of isometries between unit spheres]{On the extension of isometries between the unit spheres of a C$^*$-algebra and $B(H)$}

\author[F.J. Fern\'{a}ndez-Polo]{Francisco J. Fern\'{a}ndez-Polo}
\author[A.M. Peralta]{Antonio M. Peralta}

\address{Departamento de An{\'a}lisis Matem{\'a}tico, Facultad de
Ciencias, Universidad de Granada, 18071 Granada, Spain.}
\email{pacopolo@ugr.es}
\email{aperalta@ugr.es}




\date{December 2016}

\begin{abstract} Given two complex Hilbert spaces $H$ and $K$, let $S(B(H))$ and $S(B(K))$ denote the unit spheres of the C$^*$-algebras $B(H)$ and $B(K)$ of all bounded linear operators on $H$ and $K$, respectively. We prove that every surjective isometry $f: S(B(K)) \to S(B(H))$ admits an extension to a surjective complex linear or conjugate linear isometry $T: B(K)\to B(H)$. This provides a positive answer to Tingley's problem in the setting of $B(H)$ spaces.
\end{abstract}

\maketitle
\thispagestyle{empty}

\section{Introduction}

Let $X$ and $Y$ be normed spaces, whose unit spheres are denoted by $S(X)$ and $S(Y)$, respectively. Suppose $T: X\to Y$ is a surjective real linear isometry. The restriction $T|_{S(X)} : S(X) \to S(Y)$ defines a surjective isometry. The so-called \emph{Tingley's problem}, named after the contribution of D. Tingley \cite{Ting1987}, asks if every surjective isometry $f:S(X)\to S(Y)$ arises in this way, or equivalently, if every surjective isometry $f : S(X) \to S(Y )$ admits an extension to a surjective real linear isometry $T: X\to Y$. Tingley's achievements show that, for finite dimensional normed spaces $X$ and $Y$, every surjective isometry $f: S(X)\to S(Y)$ satisfies $f(-x) = -f(x)$ for every $x\in S(X)$ (see \cite[THEOREM in page 377]{Ting1987}).\smallskip

A solution to Tingley´s problem has been pursued by many researchers since 1987. Positive answers to Tingley's problem have been established for $\ell^p (\Gamma)$ spaces with $1\leq p\leq \infty$ (see \cite{Ding2002,Di:p,Di:8} and \cite{Di:1}), $L^{p}(\Omega, \Sigma, \mu)$ spaces, where $(\Omega, \Sigma, \mu)$ is a $\sigma$-finite measure space and $1\leq p\leq \infty$ (compare \cite{Ta:8, Ta:1} and \cite{Ta:p}), and $C_0(L)$ spaces (see \cite{Wang}). Tingley's problem also admits a positive solution in the case of finite dimensional polyhedral Banach spaces (see \cite{KadMar2012}). The reader is referred to the surveys \cite{Ding2009} and \cite{YangZhao2014} for additional details.\smallskip

In the non-commutative setting, Tingley's problem has been solved for surjective isometries between the unit spheres of two finite dimensional C$^*$-algebras (see \cite{Tan2016-2}) and for surjective isometries between the unit spheres of two finite von Neumann algebras \cite{Tan2016preprint}. A more recent contribution solves Tingley's problem for surjective isometries between the unit spheres of spaces, $K(H),$ of compact linear operators on a complex Hilbert space $H$, or more generally, for surjective isometries between the unit spheres of two compact C$^*$-algebras \cite[Theorem 3.14]{PeTan16}. The novelties in \cite{PeTan16} are based on the use of techniques of JB$^*$-triples, and Tingley's problem is also solved for surjective isometries between the unit spheres of two weakly compact JB$^*$-triples of rank greater than or equal to 5. In \cite{FerPe17} we establish a complete solution to Tingley's problem for arbitrary weakly compact JB$^*$-triples.\smallskip

Tingley's problem for surjective isometries between the unit spheres of two $B(H)$ spaces seems to be the last frontier in the studies on Tingley's problem. This paper is devoted to provided a complete solution in this case.\smallskip

The results in \cite{Tan2016-2,PeTan16,FerPe17} are based, among other techniques, on the results describing the (maximal) norm closed proper faces of the closed unit ball of a C$^*$-algebra (see \cite{AkPed92}) or of a JB$^*$-triple (see \cite{EdFerHosPe2010}).
Throughout the paper, the closed unit ball of a normed space $X$ will be denoted by $\mathcal{B}_X$. It is shown in \cite{Tan2016-2,PeTan16,FerPe17} that for a compact C$^*$-algebra $A$ (respectively, a weakly compact JB$^*$-triple $E$) the norm closed faces of $\mathcal{B}_{A}$ are determined by finite rank partial isometries in $A$ (respectively, by finite rank tripotents in $E$). However, for a general C$^*$-algebra $A$ the maximal proper faces of $\mathcal{B}_{A}$ are determined by minimal partial isometries in $A^{**}$ (see Section \ref{Sec:2} for more details). This is a serious obstacle which makes invalid the arguments in \cite{PeTan16,FerPe17} in the case of $B(H)$.\smallskip

To avoid the difficulties mentioned in the previous paragraph, our first geometric result shows that a surjective isometry $f$ from the unit sphere of a C$^*$-algebra $A$ onto the unit sphere of $B(H)$ maps minimal partial isometries in $A$ into minimal partial isometries in $B(H)$ (see Theorem \ref{t surjective isometries map minimal partial isometries into minimal partial isometries}). Apart from the just commented geometric tools, our arguments are based on techniques of functional analysis and linear algebra. In our main result we prove that given two complex Hilbert spaces $H$ and $K$, every surjective isometry $f: S(B(K)) \to S(B(H))$ admits an extension to a surjective complex linear or conjugate linear isometry $T: B(K) \to B(H)$ (see Theorem \ref{t Tingley BH spaces}). In the final result we show that the same conclusion remains true when $B(H)$ spaces are replaced by $\ell_{\infty}$-sums of $B(H)$ spaces (see Theorem \ref{thm Tyngley ellinfty sums}). The next natural question beyond these conclusions is whether Tingley's problem admits or not a positive answer for Cartan factors and atomic JBW$^*$-triples.\smallskip

It should be remarked here that the solution to Tingley's problem for surjective isometries between the unit spheres of $K(H)$-spaces in \cite{PeTan16,FerPe17} and the solution presented in this note for surjective isometries between the unit spheres of $B(H)$-spaces are completely independent results.  


\section{Surjective isometries between the unit spheres of two C$^*$-algebras}\label{Sec:2}

In this section we carry out an study of the geometric properties of surjective isometries between the unit spheres of two C$^*$-algebras with special interest on C$^*$-algebras of the form $B(H)$. We begin by gathering some technical results and concepts needed for later purposes.

\begin{proposition}\label{p preserves maximal convex subsets}{\rm(}\cite[Lemma 5.1]{ChengDong} and \cite[Lemma 3.5]{Tan2014}{\rm)} Let $X$, $Y$ be Banach spaces, and let $T : S(X) \to S(Y )$ be a surjective isometry. Then $C$ is a
maximal convex subset of $S(X)$ if and only if $T(C)$ is that of $S(Y)$. Then $C$ is a maximal proper {\rm(}norm closed{\rm)} face of $\mathcal{B}_X$ if and only if $f(C)$ is a maximal proper {\rm(}norm closed{\rm)} face of $\mathcal{B}_Y$.$\hfill\Box$
\end{proposition}

An interesting generalization of the Mazur-Ulam theorem was established by P. Mankiewicz in \cite{Mank1972}, who proved that, given two convex bodies $V\subset X$ and $W\subset Y$, every surjective isometry $g$ from $V$ onto $W$ can be uniquely extended to an affine isometry from $X$ onto $Y$. Consequently, every surjective isometry between the closed unit balls of two Banach spaces $X$ and $Y$ extends uniquely to a real linear isometric isomorphism from $X$ into $Y$.\smallskip

Let $a$ and $b$ be two elements in a C$^*$-algebra $A$. We recall that $a$ and $b$ are orthogonal ($a\perp b$ in short) if $ ab^* = b^* a =0$. Symmetric elements in $A$ are orthogonal if and only if their product is zero.\smallskip

For each element $a$ in a C$^*$-algebra $A$, the symbol $|a|$ will denote the element $(a^* a)^{\frac12}\in A$. Throughout this note, for each $x\in A$, $\sigma (x)$ will denote the spectrum of the element $x$. We observe that $\sigma(|a|)\cup \{0\} = \sigma(|a^*|)\cup \{0\}$, for every $a\in A$. Let $a = v |a|$ be the polar decomposition of $a$ in $A^{**}$, where $v$ is a partial isometry in $A^{**}$, which, in general, does not belong to $A$ (compare \cite{S}). It is further known that $v^*v$ is the range projection of $|a|$ ($r(|a|)$ in short), and for each $h\in C(\sigma(|a|)),$ with $h(0)=0$ the element $v h(|a|)\in A$ (see \cite[Lemma 2.1]{AkPed77}).\smallskip

Proposition \ref{p preserves maximal convex subsets} points out the importance of an appropriate description of the maximal proper faces of the closed unit ball $\mathcal{B}_A$ of a C$^*$-algebra $A$. A complete study was established by C.A. Akemann and G.K. Pedersen in \cite{AkPed92}. When $A$ is a von Neumann algebra, weak$^*$-closed faces in $\mathcal{B}_A$ were originally determined by C.M. Edwards and G.T. R\"{u}ttimann in \cite{EdRutt88}, who proved that general weak$^*$-closed faces in $\mathcal{B}_A$ have the form $$F_v=v + (1 - vv^*) \mathcal{B}_{A} (1 - v^*v) = \{x\in \mathcal{B}_A:\ xv^* = vv^*\},$$ for some partial isometry $v$ in $A$. Actually, the mapping $v\mapsto F_v$ is an anti-order isomorphism from the complete lattice of partial isometries in $A$ onto the complete lattice of weak$^*$-closed faces of $\mathcal{B}_A$, where the partial order in the set of partial isometries of $A$ is given by $v\leq u$ if and only if $u = v + (1-vv^*) u (1-v^*v)$ (see \cite[Theorem 4.6]{EdRutt88}).\smallskip

However, partial isometries in a general C$^*$-algebra $A$ are not enough to determine all the norm-closed faces in $\mathcal{B}_A$, even more after recalling the existence of C$^*$-algebras containing no partial isometries. In the general case, certain partial isometries in the second dual $A^{**}$ are required to determine the facial structure of $\mathcal{B}_A$. We recall that a projection $p$ in $A^{**}$ is called \emph{open} if $A\cap (p A^{**} p)$ is weak$^*$-dense in $p A^{**} p$ (see \cite[\S 3.11]{Ped} and \cite[\S III.6]{Tak}). A projection $p \in A^{**}$ is said to be \emph{closed} if $1-p$ is open. A closed projection $p\in A^{**}$ is \emph{compact} if $p\leq x$ for some positive norm-one element $x \in A$. A partial isometry $v\in A^{**}$ \emph{belongs locally to $A$} if $v^*v$ is a compact projection and there exists a norm-one element $x$ in $A$ satisfying $v = x v^*v$ (compare \cite[Remark 4.7]{AkPed92}). 
\smallskip

It is shown in \cite[Lemma 4.8 and Remark 4.11]{AkPed92} that ``the partial isometries that belong locally to $A$ are obtained by taking an element $x$ in $A$ with norm 1 and polar decomposition $x = u |x|$ (in $A^{**}$), and then letting $v = ue$ for some compact projection $e$ contained in the spectral projection $\chi_{_{\{1\}}}(|x|)$ of $|x|$ corresponding to the eigenvalue 1.''\smallskip

It should be noted that a partial isometry $v$ in $A^{**}$ belongs locally to $A$ if and only if it is compact in the sense introduced by C.M. Edwards and G.T. R\"{u}ttimann in \cite[Theorem 5.1]{EdRu96}.\smallskip

The facial structure of the unit ball of a C$^*$-algebra is completely described by the following result due to C.A. Akemann and G.K. Pedersen.

\begin{theorem}\label{t faces AkPed}\cite[Theorem 4.10]{AkPed92} Let $A$ be a C$^*$-algebra. The norm closed faces of the unit ball of $A$ have the form $$F_v=\left(v + (1 - vv^*) \mathcal{B}_{A^{**}} (1 - v^*v)\right)\cap \mathcal{B}_{A} = \{x\in \mathcal{B}_A:\ xv^* = vv^*\},$$ for some partial isometry $v$ in $A^{**}$ belonging locally to $A$. Actually, the mapping $v\mapsto F_v$ is an anti-order isomorphism from the complete lattice of partial isometries in $A^{**}$ belonging locally to $A$ onto the complete lattice of norm closed faces of $\mathcal{B}_A$. $\hfill\Box$
\end{theorem}

A non-zero partial isometry $e$ in a C$^*$-algebra $A$ is called minimal if $ee^*$ (equivalently, $e^* e$) is a minimal projection in $A$, that is, $ ee^* A e e^* =\mathbb{C}  ee^*.$ By Kadison's transitivity theorem minimal partial isometries in $A^{**}$ belong locally to $A$, and hence every maximal proper face of the unit ball of a C$^*$-algebra $A$ is of the form $$\left(v + (1 - vv^*) \mathcal{B}_{A^{**}} (1 - v^*v)\right)\cap \mathcal{B}_{A}$$ for a unique minimal partial isometry $v$ in $A^{**}$ (compare \cite[Remark 5.4 and Corollary 5.5]{AkPed92}).\smallskip

Our main goal in this section is to show that a surjective isometry $f : S(A)\to S(B)$ between the unit spheres of two C$^*$-algebras maps minimal partial isometries into minimal partial isometries. In a first step we shall show that, for each minimal partial isometry $e$ in $A$, 1 is isolated in the spectrum of $|f(e)|$.

\begin{theorem}\label{t surjective isometries map minimal partial isometries into points of strong subdiff}
Let $A$ and $B$ be C$^*$-algebras, and suppose that $f: S(A) \to S(B)$ is a surjective isometry. Let $e$ be a minimal partial isometry in $A$. Then $1$ is isolated in the spectrum of $|f(e)|$.
\end{theorem}

\begin{proof} Since $e$ also is a minimal partial isometry in $A^{**}$ and belongs (locally) to $A$,  the set $F_{e} =e + (1 - ee^*) \mathcal{B}_A (1 - e^*e)$ is a maximal proper face of $\mathcal{B}_{A}$. Applying Proposition \ref{p preserves maximal convex subsets} and Theorem \ref{t faces AkPed} we deduce the existence of a minimal partial isometry $w$ in $B^{**}$ such that \begin{equation}\label{eq faces thm1} f(F_e) =\! F_{w} =\!\left(w + (1 - ww^*) \mathcal{B}_{B^{**}} (1 - w^*w)\right)\cap \mathcal{B}_B \!=\{x\in \mathcal{B}_A: xw^* = ww^*\}.
\end{equation} Since $f(e)\in f(F_e) = F_{w}$ we have $f(e) = w + (1 - ww^*) f(e) (1 - w^*w)$.\smallskip

Arguing by contradiction, we assume that $1$ is not isolated in $\sigma( |f(e)|)$. Let $f(e) = r |f(e)|$ denote the polar decomposition of $f(e)$.\smallskip

By assumptions we can find $t_0\in \sigma(|f(e)|)$ satisfying $\frac{3}{\sqrt{10}}<{t_0}<1$. Let us consider the functions ${h_1}$ and $h_2$ in the unit sphere of $C_0(\sigma(|f(e)|))$ given by
$$h_1(t):=\left\{%
\begin{array}{ll}
    \frac{t}{t_0}, & \hbox{if $0\leq t\leq t_0$} \\
    \hbox{affine}, & \hbox{if $t_0\leq t\leq \frac 12 (t_0+1)$} \\
    0, & \hbox{if $\frac 12 (t_0+1)\leq t\leq 1$} \\
\end{array}%
\right. \  ;\ \  h_2(t):=\left\{%
\begin{array}{ll}
    0, & \hbox{if $0\leq t\leq \frac 12 (t_0+1)$} \\
    \hbox{affine}, & \hbox{if $\frac 12 (t_0+1)\leq t\leq 1$} \\
    1, & \hbox{if $ t=1$.} \\
\end{array}%
\right.$$ We set $\widehat{x} = r h_1 (|f(e)|)$ and $\widehat{y} = r h_2 (|f(e)|)$. Obviously $h_1 (|f(e)|)$ and $h_2 (|f(e)|)$ are positive elements in $S(B)$ satisfying $h_1 (|f(e)|) h_2 (|f(e)|)=0.$ Since $$\widehat{x} \widehat{y}^* =  r h_1 (|f(e)|) h_2 (|f(e)|) r^* =0,$$ and $$ \widehat{y}^* \widehat{x} =  h_2 (|f(e)|) r^* r h_1 (|f(e)|)  =h_2 (|f(e)|) h_1 (|f(e)|) =0,$$ it follows that $\hat{x}\perp\hat{y}.$\smallskip

Let $x=f^{-1}(-\hat{x})\in S(A)$ and $y=f^{-1}(\hat{y})\in S(A)$. Since $f$ is an isometry we deduce that $$1=\|\hat{x}+\hat{y}\|=\|\hat{y}- (-\hat{x})\|= \|f(y)-f(x)\|= \|y-x\|,$$ and $$1+t_0=\|f(e)+\hat{x}\|=\|f(e)-f(x)\|=\|e-x\|.$$\smallskip

We recall that, from \eqref{eq faces thm1}, $f(e) = w + k$ where $k=(1 - ww^*) f(e) (1 - w^*w)$ satisfies $k^* w = 0 = w k^*$, which proves that $k\perp w$. Let $r_0$ denote the (unique) partial isometry appearing in the polar decomposition of $k$. Since $r$ is the partial isometry in the polar decomposition of $f(e)$, $w\perp k$, and $f(e) = w + k$, it follows that $r = w + r_0$ with $r_0\perp w$. We also know that\label{page ref face} $|f(e)| = w^*w + |k|$, and hence a simple application of the continuous functional calculus (having in mind that $h_2 (1) =1$) shows that $h_2 (|f(e)|) = w^* w + h_2 (|k|)$, with $w^* w \perp h_2 (|k|)$. We therefore have \begin{equation}\label{eq 1 030117} \hat{y} w^* = r h_2 (|f(e)|) w^* = r w^* w w^*+ r h_2 (|k|) w^* = r w^* =(w+r_0) w^* = w w^*,
 \end{equation} which implies that $\hat{y} \in F_w,$ and consequently $y =f^{-1}(\hat{y})\in F_e$ (see \eqref{eq faces thm1}).\smallskip

We claim that \begin{equation}\label{eq claim in thm 1} \|  ee^* (e-x) e^*e\|>1.
\end{equation}

The element $e-x$ has norm $1+t_0>1$. Suppose that $\{H_i\}_{I}$ is a family of complex Hilbert spaces and $\pi : A\to \bigoplus_i^{\ell_{\infty}} B(H_i)$ is an isometric $^*$-homomorphism with weak$^*$-dense range (we can consider, for example, the \emph{atomic representation} of $A$ \cite[4.3.7]{Ped}, where the family $I$ is precisely the set of all pure states of $A$ and $\pi$ is the direct sum of all the irreducible representations associated with the pure states \cite[Theorem 3.13.2]{Ped}). For each $j\in I$, let $P_j$ denote the projection of $\bigoplus_i^{\ell_{\infty}} B(H_i)$ onto $B(H_j)$ and let $\pi_j = P_j \circ \pi$. Clearly, $\pi_j$ is a $^*$-homomorphism with weak$^*$-dense range. Since $e$ is a minimal partial isometry, there exists a unique $i_0\in I$ such that $\pi_{i_0} (e)$ is a non-zero (minimal) partial isometry and $\pi_j (e) =0$, for every $j\neq i_0$. We also know that $\|x\|=1$, and thus $\| \pi_{i_0} (e-x)\|= 1+{t_0}.$\smallskip

Let $\pi_{i_0} (e-x)= u |\pi_{i_0} (e-x)|$ be the polar decomposition of $\pi_{i_0} (e-x)$ in $B(H_{i_0})$. Take $0<\varepsilon <t_0 -\frac{3}{\sqrt{10}}$.  Since $\| |\pi_{i_0} (e-x)| \|=1+t_0$ we can find a minimal projection $q= \xi\otimes \xi\in B(H_{i_0})$ with $\|\xi\|=1$ in $H_{i_0}$ satisfying $q\leq u^* u$ and \begin{equation}\label{eq eval at xi big} 1+{t_0} -\varepsilon <\langle |\pi_{i_0} (e-x)| (\xi) / \xi\rangle ,
\end{equation} and \begin{equation}\label{eq eval at xi big 2}  1+{t_0} -\varepsilon < \langle |\pi_{i_0} (e-x)| (\xi) / \xi\rangle \leq \langle |\pi_{i_0} (e-x)| (\xi) /  |\pi_{i_0} (e-x)| (\xi)\rangle^{\frac12} \ \|\xi\|
\end{equation} $$= \langle |\pi_{i_0} (e-x)|^2 (\xi) /  \xi \rangle^{\frac12} = \langle \pi_{i_0} (e-x)^* \pi_{i_0} (e-x) (\xi) /  \xi \rangle^{\frac12}.$$ The element $v = u q$ is a minimal partial isometry in $B(H_{i_0})$.\smallskip

We observe that $\pi (e) = \pi_{i_0} (e)$ and $v$ are not orthogonal. Otherwise, $\pi_{i_0} (e)^* \pi_{i_0} (e) \perp v^* v=q$, and hence $\pi_{i_0} (e) q = 0 = q \pi_{i_0} (e)^*,$ which, by \eqref{eq eval at xi big 2}, implies that
$$\left(1+{t_0} -\varepsilon\right)^2< \langle\pi_{i_0} (e-x)^* \pi_{i_0} (e-x) (\xi) / \xi\rangle = \langle q \pi_{i_0} (e-x)^* \pi_{i_0} (e-x) q (\xi) / \xi\rangle$$
$$ = \langle q \pi_{i_0} (x)^* \pi_{i_0} (x) q (\xi) / \xi\rangle \leq \| \pi_{i_0} (x)^* \pi_{i_0} (x) \| = \|\pi_{i_0} (x)\|^2 \leq \|x\|^2 = 1,$$ which is impossible.\smallskip

Therefore, $\pi_{i_0} (e)$ and $v$ are two minimal partial isometries in $B(H_{i_0})$ which are not orthogonal. They must be of the form $\pi_{i_0} (e) = \eta_1 \otimes \xi_1$ and $v= \widetilde{\eta}_1 \otimes \widetilde{\xi}_1$ for suitable $\xi_1,\eta_1, \widetilde{\xi}_1,\widetilde{\eta}_1\in S(H_{i_0})$ with $|\langle \xi_1 / \widetilde{\xi}_1 \rangle|+ |\langle  \widetilde{\eta}_1 / \eta_1 \rangle|\neq 0$. Let us consider two orthonormal systems $\{\eta_1, \eta_2\}$ and $\{\xi_1, \xi_2\}$ such that
$$v=\alpha \pi_{i_0} (e)+ \beta v_{12}+ \delta v_{22}+\gamma v_{21},$$ where $\pi_{i_0} (e):= v_{11}= \eta_1 \otimes \xi_1$, $v_{12} = {\eta}_2 \otimes {\xi}_1$, $v_{21} = {\eta}_1 \otimes {\xi}_2$, $v_{22} = {\eta}_2 \otimes {\xi}_2$,
$\alpha= \langle \xi_1 / \widetilde{\xi}_1 \rangle \langle  \widetilde{\eta}_1 / \eta_1 \rangle,$ $ \beta = \langle \xi_1 / \widetilde{\xi}_1 \rangle \langle  \widetilde{\eta}_1 / \eta_2 \rangle,$  $\gamma= \langle \xi_2 / \widetilde{\xi}_1 \rangle \langle  \widetilde{\eta}_1 / \eta_1 \rangle,$ $\delta = \langle \xi_2 / \widetilde{\xi}_1 \rangle \langle  \widetilde{\eta}_1 / \eta_2 \rangle \in \mathbb{C}.$ It is easy to check that $|\alpha|^2+ |\beta|^2+|\gamma|^2+ |\delta|^2$ $=|\langle \xi_1 / \widetilde{\xi}_1 \rangle|^2 \|\widetilde{\eta}_1\|^2+ |\langle \xi_2 / \widetilde{\xi}_1 \rangle|^2 \|\widetilde{\eta}_1 \|^2 = \| \widetilde{\xi}_1\|^2=1$, and $\alpha \delta= \beta \gamma$.\label{eq min partial isometry as a matrix}\smallskip

For each $(i,j)\in \{1,2\}^2$, let $\varphi_{ij}\in S(B(H_{i_0})_*)$ be the unique extreme point of the unit ball $\mathcal{B}_{B(H_{i_0})_*}$ defined by $\varphi_{ij} (z) := \langle z(\xi_i) / \eta_j\rangle$ ($z\in B(H_{i_0})$). We shall also consider $\varphi_v \in S(B(H_{i_0})_*)$, defined by $\varphi_{v} (z) := \langle z(\widetilde{\xi}_1) / \widetilde{\eta}_1\rangle$ ($z\in B(H_{i_0})$).  Each $\varphi_{ij}$ is supported by $v_{ij}\in S(B(H_{i_0}))$, while $\varphi_v$ is supported by $v$.\smallskip

Clearly,  the identity $$\varphi_{v}(\pi_{i_0} (e))= \langle \pi_{i_0} (e)(\widetilde{\xi}_1) / \widetilde{\eta}_1\rangle = \langle \widetilde{\xi}_1 / \xi_1 \rangle \langle \eta_1  / \widetilde{\eta}_1 \rangle =\bar{\alpha}$$ holds, and similarly we have $$\varphi_{v}\pi_{i_0}(x)=\bar{\alpha} z_{11}+ \bar{\beta} z_{12}+ \bar{\delta} z_{22}+\bar{\gamma} z_{21},$$ where $z_{ij}=\varphi_{ij} \pi_{i_0}(x)$ for all $(i,j)\in \{1,2\}^2$. We also know that $\widetilde{\xi}_1\otimes \widetilde{\xi}_1= v^* v =q =\xi\otimes \xi \leq u^* u,$ and thus $ \widetilde{\xi}_1 = \mu_0 \xi$, for a suitable $\mu_0\in \mathbb{C}$ with $|\mu_0|=1$. We deduce from \eqref{eq eval at xi big} that $$1+{t_0} -\varepsilon <\langle |\pi_{i_0} (e-x)| (\xi) / \xi\rangle = \langle u^*u |\pi_{i_0} (e-x)| (\xi) / \xi\rangle = \langle u |\pi_{i_0} (e-x)| (\xi) / u q (\xi) \rangle$$ $$= \langle \pi_{i_0} (e-x) (\xi) / v (\xi) \rangle = \langle \pi_{i_0} (e-x) (\widetilde{\xi}_1) / v (\widetilde{\xi}_1) \rangle = \langle \pi_{i_0} (e-x) (\widetilde{\xi}_1) / \widetilde{\eta}_1 \rangle =\varphi_v (\pi_{i_0} (e-x))$$ $$  =\varphi_v (\pi_{i_0} (e)) + \varphi_v (\pi_{i_0} (x)) \leq |\alpha| + |\varphi_v (\pi_{i_0} (x)) |\leq |\alpha| +1, $$ which proves $$ {t_0} -\varepsilon< |\alpha|\leq 1.$$\smallskip

Now, the equality $|\alpha|^2+ |\beta|^2+|\gamma|^2+ |\delta|^2=1$ implies that $$|\beta|^2,|\gamma|^2, |\delta|^2\leq 1-({t_0} -\varepsilon)^2<\frac{1}{10},$$ and since $|z_{ij}|\leq 1$, we have $$1+{t_0} -\varepsilon< \varphi_v \pi_{i_0} (e-x) = |\varphi_v \pi_{i_0} (e-x)|=|\bar{\alpha}-( \bar{\alpha} z_{11}+ \bar{\beta} z_{12}+ \bar{\delta} z_{22}+\bar{\gamma} z_{21}) | $$ $$\leq |\alpha| |1-z_{11}|+ |z_{12}||\beta| + |z_{21}||\gamma|+ |z_{22}||\delta| \leq |1-z_{11}|+ |\beta| + |\gamma|+ |\delta|\leq |1-z_{11}|+\frac{3}{\sqrt{10}}.$$

Let us observe that $\pi_{i_0} (e) \pi_{i_0} (e)^* = v_{11}v_{11}^*= \eta_1 \otimes \eta_1$ and  $ \pi_{i_0} (e)^*  \pi_{i_0} (e)= v_{11}^*v_{11}= \xi_1 \otimes \xi_1$. Therefore $$1<1+{t_0} -\varepsilon -\frac{3}{\sqrt{10}}\leq |1-z_{11}| = |\varphi_{11} \pi_{i_0} (e-x) |$$ $$ = \Big|\varphi_{11} \Big((\pi_{i_0} (e) \pi_{i_0} (e)^*) \pi_{i_0} (e-x) (\pi_{i_0} (e)^* \pi_{i_0} (e)) \Big)\Big|$$ $$\leq \Big\| (\pi_{i_0} (e) \pi_{i_0} (e)^*) \pi_{i_0} (e-x) (\pi_{i_0} (e)^* \pi_{i_0} (e)) \Big\| $$ $$=  \| \pi_{i_0} ( (e e^*) (e-x) (e^* e)) \|\leq \| (e e^*) (e-x) (e^* e)\|,$$ which proves the claim in \eqref{eq claim in thm 1}.\smallskip

Finally, since $y\in F_e = e + (1 - ee^*) \mathcal{B}_A (1 - e^*e)$ we can write $$y = e + (1 - ee^*) y (1 - e^*e),$$ and we deduce from \eqref{eq claim in thm 1} that $$1 = \|y-x\| \geq \| ee^* (y-x) e^* e \| = \| ee^* (e + (1 - ee^*) y (1 - e^*e)-x) e^* e \| $$ $$= \| ee^* (e + x) e^* e \| >1,$$ leading to the desired contradiction.\end{proof}

The problem of dealing with minimal faces of the unit ball of a C$^*$-algebra $A$ is that we need to handle minimal partial isometries in $A^{**}$ (compare Theorem \ref{t faces AkPed}). We present now a technical result which will be used later to facilitate the arguments depending on the facial structure of $\mathcal{B}_A$.

\begin{lemma}\label{l minimal in A** not in A}
Let $A$ be a C$^*$-algebra. The following statements hold:\begin{enumerate}[$(a)$]\item Every minimal projection $p$ in $A^{**}\backslash A$ is orthogonal to all minimal projections in $A$;
\item Every minimal partial isometry $u$ in $A^{**}\backslash A$ is orthogonal to all minimal partial isometries in $A$.
\end{enumerate}
\end{lemma}

\begin{proof} $(a)$ Suppose $p$ is a minimal projection in $A^{**}\backslash A$. Let $q$ denote a minimal projection in $A$. Arguing by contradiction we assume that $p q\neq 0$.\smallskip

As in the proof of Theorem \ref{t surjective isometries map minimal partial isometries into points of strong subdiff} let $\pi : A\to \bigoplus_i^{\ell_{\infty}} B(H_i)$ be an isometric $^*$-homomorphism with weak$^*$-dense range, where $\{H_i\}_{I}$ is a family of complex Hilbert spaces (consider, for example, the \emph{atomic representation} of $A$ \cite[4.3.7]{Ped}). By the weak$^*$-density of $A$ in $\bigoplus_i^{\ell_{\infty}} B(H_i)$ and the separate weak$^*$-continuity of the product of every von Neumann algebra, $\pi (q)$ is a minimal projection in $\bigoplus_i^{\ell_{\infty}} B(H_i)$. Clearly, the images of the mappings $L_{\pi(q)}: x \mapsto \pi(q) x$ and $R_{\pi(q)}: x\mapsto x \pi(q)$ ($\forall x\in \bigoplus_i^{\ell_{\infty}} B(H_i)$) are contained in suitable Hilbert spaces. It follows that the left and right multiplication operators $L_q$ and $R_q$ by $q$ on $A$ factors through a Hilbert space, and thus they are weakly compact (compare \cite{DaFiJoPe}). Consequently, the spaces $(1-q) A q$, $q A (1-q)$ and $q A q = \mathbb{C} q$ are all reflexive. Applying the Krein-\v{S}mulian theorem we deduce that $(1-q) A q$, $q A (1-q)$ and $q A q = \mathbb{C} q$ are weak$^*$-closed in $A^{**}$, showing that $$(1-q) A^{**} q = (1-q) A q,\ q A^{**} (1-q) = q A (1-q) \subseteq A, \hbox{ and }q A^{**} q = q A q  = \mathbb{C} q.$$

We recall now an useful matricial representation theorem. Let $C$ denote the C$^*$-subalgebra of $A^{**}$ generated by $p$ and $q$. Since $p$ and $q$ are minimal projections in $A^{**}$, Theorem 1.3 in \cite{RaSin} (see also \cite[\S 3]{Ped68}) assures the existence of $t\in [0,1]$ and a $^*$-isomorphism $\Phi : C \to M_2 (\CC)$ such that $\Phi (q) = \left( \begin{array}{cc}
                                                     1 & 0 \\
                                                     0 & 0 \\
                                                   \end{array}
                                                 \right)$ and  $\Phi (p) = \left(
                                                                             \begin{array}{cc}
                                                                               t & \sqrt{t(1-t)} \\
                                                                               \sqrt{t(1-t)} & 1-t \\
                                                                             \end{array}
                                                                           \right).$ Since $pq \neq 0$ we know that $t\neq 0$. Clearly, $\Phi^{-1} \left(
                                                                             \begin{array}{cc}
                                                                               0 & 1 \\
                                                                               1 & 0 \\
                                                                             \end{array}
                                                                           \right)\in q C (1-q) \oplus (1-q) C q \subset A,$ and $\Phi^{-1} \left(
                                                                             \begin{array}{cc}
                                                                               1 & 0 \\
                                                                               0 & 0 \\
                                                                             \end{array}
                                                                           \right) =q \in A$.
Then $$\Phi ^{-1} \left(
                                                                             \begin{array}{cc}
                                                                               0 & 0 \\
                                                                               0 & 1 \\
                                                                             \end{array}
                                                                           \right) = \Phi ^{-1} \left( \left(
                                                                             \begin{array}{cc}
                                                                               0 & 1 \\
                                                                               1 & 0 \\
                                                                             \end{array}
                                                                           \right) \left(
                                                                             \begin{array}{cc}
                                                                               1 & 0 \\
                                                                               0 & 0 \\
                                                                             \end{array}
                                                                           \right) \left(
                                                                             \begin{array}{cc}
                                                                               0 & 1 \\
                                                                               1 & 0 \\
                                                                             \end{array}
                                                                           \right)\right)\in \Phi^{-1} (\Phi(A\cap C)) =A\cap C.$$ By linearity $p = \Phi^{-1} \left(
                                                                             \begin{array}{cc}
                                                                               t & \sqrt{t(1-t)} \\
                                                                               \sqrt{t(1-t)} & 1-t \\
                                                                             \end{array}
                                                                           \right)\in A$ which is impossible.\smallskip

$(b)$ Suppose now that $u$ is a minimal partial isometry in $A^{**}\backslash A$ and $v$ is a minimal partial isometry in $A$. We shall first show that $uu^*, u^*u\in A^{**}\backslash A$. Indeed, since every minimal partial isometry in $A^{**}$ belongs locally to $A$ (compare Kadison's transitivity theorem and \cite[Remark 5.4 and Corollary 5.5]{AkPed92}), there exists a norm one element $x\in A$ satisfying $x= u + (1-uu^*) x (1-u^*u)$. If $uu^*$ (respectively, $u^* u$) lies in $A$ then $u = uu^* x \in A$ (respectively, $u = x u^* u \in A$) which is impossible.\smallskip

We have therefore shown that $uu^*, u^*u\in A^{**}\backslash A$ are minimal projections, while $vv^*, v^*v$ are minimal projections in $A$. It follows from $(a)$ that $uu^*, u^*u \perp vv^*, v^*v$. Finally, the identities $u^* v = u^* u u^* v v^* v =0$ and $v u^* = v v^* v u^* u u^* =0$ prove that $u\perp v$.\end{proof}

We are now in position to show that a surjective isometry between the unit spheres of two C$^*$-algebras maps minimal partial isometries to minimal partial isometries.

\begin{theorem}\label{t surjective isometries map minimal partial isometries into minimal partial isometries}
Let $A$ be a C$^*$-algebra, and let $H$ be a complex Hilbert space. Suppose that $f: S(A) \to S(B(H))$ is a surjective isometry. Let $e$ be a minimal partial isometry in $A$. Then $f(e)$ is a minimal partial isometry in $B(H)$. Moreover, there exits a surjective real linear isometry $$T_{e} : (1 - ee^*) A  (1 - e^*e)\to  {(1 - f(e)f(e)^*) B(H) (1 - f(e)^*f(e))}$$ such that $$ f(e + x) = f(e) + T_e (x), \hbox{ for all $x$ in } \mathcal{B}_{(1 - ee^*) A  (1 - e^*e)}.$$ In particular the restriction of $f$ to the face $F_{e} =e + (1 - ee^*) \mathcal{B}_A (1 - e^*e)$ is a real affine function.
\end{theorem}

\begin{proof}Arguing as in the beginning of the proof of Theorem \ref{t surjective isometries map minimal partial isometries into points of strong subdiff}, the set $$F_{e} =e + (1 - ee^*) \mathcal{B}_A (1 - e^*e)$$ is a maximal proper face of $\mathcal{B}_{A}$, and thus, by Proposition \ref{p preserves maximal convex subsets} and Theorem \ref{t faces AkPed}, there exists a minimal partial isometry $w$ in $B(H)^{**}$ such that \begin{equation}\label{eq faces thm2} f(F_e) = F_{w} =\left(w + (1 - ww^*) \mathcal{B}_{B(H)^{**}} (1 - w^*w)\right)\cap \mathcal{B}_{B(H)}.\end{equation}

\emph{We claim that $w\in B(H)$}. Suppose, on the contrary that $w\in B(H)^{**}\backslash B$.\smallskip

Theorem \ref{t surjective isometries map minimal partial isometries into points of strong subdiff} implies that 1 is an isolated point in $\sigma (|f(e)|)$, and hence the function $\chi_{_{\{1\}}}$ belongs to $C_0(\sigma (|f(e)|))$. Let $f(e) = r |f(e)|$ denote the polar decomposition of $f(e)$. An application of the continuous functional calculus proves that $\hat{v}= r \ \chi_{_{\{1\}}} (|f(e)|)$ is a partial isometry in $B(H).$ Furthermore, since $f(e)\in f(F_e) = F_w$, we deduce that $\hat{v} \in F_w$ and \begin{equation}\label{eq first decompostion of hatv in thm2} \hat{v} = {w} + (1-{w} {w}^*) \hat{v} (1-{w}^* {w})
 \end{equation} (compare the arguments in the proof of \eqref{eq 1 030117} in page \pageref{page ref face}).\smallskip

In $B(H)$ we can always find a minimal partial isometry $\hat{w}\in B(H)$ satisfying \begin{equation}\label{eq hatv and minimal hatw} \hat{v} = \hat{w} + (1-\hat{w} \hat{w}^*) \hat{v} (1-\hat{w}^* \hat{w}).
 \end{equation}Since, by assumptions $w\in B(H)^{**}\backslash B(H)$, Lemma \ref{l minimal in A** not in A} implies that $w\perp \hat{w}$, $$\hat{w} = (1-{w} {w}^*) \hat{w} (1-{w}^* {w})$$ and hence, by \eqref{eq first decompostion of hatv in thm2} we get $$\hat{v} - \hat{w} = {w} + (1-{w} {w}^*) (\hat{v}- \hat{w}) (1-{w}^* {w}).$$

By hypothesis, $$2=\|f(e)+\hat{w}\|=\|f(e)-(-\hat{w})\| =\|e-f^{-1}(-\hat{w})\|$$ where $e$ is a minimal partial isometry in $A$. Proposition 2.2 in \cite{FerPe17} proves that $$\hat{e}=f^{-1}(-\hat{w})=-e+(1-ee^*) \hat{e} (1-e^*e).$$

By construction $f(e) = \hat{v} + (1-\hat{v} \hat{v}^*) f(e) (1-\hat{v}^* \hat{v})$, and by \eqref{eq hatv and minimal hatw}, $$f(e) = \hat{w} + (1-\hat{w} \hat{w}^*) \hat{v} (1-\hat{w}^* \hat{w}) + (1-\hat{v} \hat{v}^*) f(e) (1-\hat{v}^* \hat{v}),$$ and consequently $\|f(e)-\hat{w}\|\leq 1$. Having in mind that $$f(e)-\hat{w} = f(e) - (1-{w} {w}^*) \hat{w} (1-{w}^* {w})\in F_w + (1-{w} {w}^*) \mathcal{B}_{B(H)} (1-{w}^* {w})$$ we get $$f(e)-\hat{w} = w + (1-w w^*) (f(e)-\hat{w}) (1-w^* w) \in F_{w}.$$ We deduce from \eqref{eq faces thm2} that $z=f^{-1}(f(e)-\hat{w})\in F_{e}$, and thus $z=e+(1-ee^*)z (1-e^*e),$ which leads to $$1=\|f(e)\|=\|f(e)-\hat{w}+\hat{w}\|=\|(f(e)-\hat{w})-(-\hat{w})\|=\|f(z) -f(\hat e)\|$$ $$=\|z-\hat{e}\|=\|e+(1-ee^*)z (1-e^*e) -(-e+(1-ee^*) \hat{e} (1-e^*e))\|$$ $$= \|2e+(1-ee^*) (z-\hat{e}) (1-e^*e) \|= 2,$$ and hence to a contradiction. Therefore $w\in B(H)$ and $$F_{w} = w + (1 - ww^*) \mathcal{B}_{B(H)} (1 - w^*w).$$

We can argue now as in the proof of \cite[Proposition 3.1]{PeTan16} to conclude. We insert a short argument here for completeness reasons. We have established that $$f\left( e +  \mathcal{B}_{(1 - ee^*) A  (1 - e^*e)} \right)=f(F_e) =F_w = w +  \mathcal{B}_{(1 - ww^*) B(H) (1 - w^*w)}.$$ Let $\mathcal{T}_{x_0}$ denote the translation with respect to $x_0$, that is $\mathcal{T}_{x_0} (x) = x+x_0$. The mapping $f_{e} = \mathcal{T}_{w}^{-1}|_{f(F_e)} \circ f|_{F_e} \circ \mathcal{T}_{e}|_{\mathcal{B}_{(1 - ee^*) A  (1 - e^*e)}}$ is a surjective isometry from $\mathcal{B}_{(1 - ee^*) A  (1 - e^*e)}$ onto $\mathcal{B}_{(1 - ww^*) B(H) (1 - w^*w)}$. Mankiewicz's theorem (see \cite{Mank1972}) implies the existence of a surjective real linear isometry $T_{e} : (1 - ee^*) A  (1 - e^*e)\to  {(1 - ww^*) B(H) (1 - w^*w)}$ such that $f_{e} = T_{e}|_{S((1 - ee^*) A  (1 - e^*e))}$ and hence $$ f(e + x) = w + T_e (x), \hbox{ for all $x$ in } \mathcal{B}_{(1 - ee^*) A  (1 - e^*e)}.$$ In particular $f(e) = w$.\smallskip

For the final statement we simple write $$f|_{F_e}= \mathcal{T}_{w}|_{\mathcal{B}_{(1 - ww^*) B(H) (1 - w^*w)}} \circ f_{e} \circ \mathcal{T}_{e}^{-1}|_{F_e}= \mathcal{T}_{w}|_{\mathcal{B}_{(1 - ww^*) B(H) (1 - w^*w)}} \circ T_{e} \circ \mathcal{T}_{e}^{-1}|_{F_e}$$ as a composition of real affine functions.\end{proof}

The next technical lemma is obtained with basic techniques of linear algebra.

\begin{lemma}\label{l minimal partial isometries} Let $H$ be a complex Hilbert space, and let $v_1,v_2, v_3$, $e_1$ and $e_2$ be minimal partial isometries in $B(H)$ satisfying $e_1\perp e_2$, $v_1\perp v_2,$ $v_1\perp v_3$, $$-v_1+v_3 = e_1 -e_2, \hbox{ and } v_1+v_2 = e_1 +e_2.$$ Then $v_1 = e_2$ and $v_2=e_1= v_3$.
\end{lemma}

\begin{proof} Since $v_1\perp v_2,v_3$, by multiplying the identities $-v_1+v_3 = e_1 -e_2$ and $v_1+v_2 = e_1 +e_2$ on the left by $v_1^*$ we get $$-v_1^* v_1 = v_1^* e_1 - v_1^* e_2,\hbox{ and } v_1^* v_1  = v_1^* e_1 + v_1^*e_2,$$ which shows that $v_1^* e_1 =0$. Multiplying by $v_1^*$ on the right we prove $e_1 v_1^* =0$. We have therefore shown that $v_1\perp e_1$.\smallskip

Applying that $e_1\perp e_2$ and $v_1\perp v_2$ we get $e_1 e_1^* + e_2 e_2^* = v_1 v_1^* + v_2 v_2^*$ where $e_1 e_1^*$ and  $v_1 v_1^*$ are orthogonal rank one projections. Thus, $e_1 e_1^* = e_1 e_1^* v_2 v_2^*$, and by minimality $v_2 v_2^* = e_1 e_1^*$. We can similarly prove $v_2^* v_2 = e_1^* e_1.$ Finally $$ v_2 = v_2 v_2^* v_2 = v_2 v_2^* (v_1+v_2) = v_2 v_2^* (e_1+e_2) = e_1 e_1^* (e_1+e_2) = e_1,$$ and the rest is clear.
\end{proof}

Next, we shall establish several consequences of the above theorem.

\begin{theorem}\label{t surjective isometries between spheres of B(H)} Let $f: S(B(K)) \to S(B(H))$ be a surjective isometry where $H$ and $K$ are complex Hilbert spaces with dimension greater than or equal to 3. Then the following statements hold:\begin{enumerate}[$(a)$]\item For each minimal partial isometry $v$ in $B(K)$, the mapping $$T_v : (1 - vv^*) B(K)  (1 - vv^*)\to  {(1 - f(v)f(v)^*) B(H) (1 - f(v)^*f(v))}$$ given by Theorem \ref{t surjective isometries map minimal partial isometries into minimal partial isometries} is complex linear or conjugate linear;
\item For each minimal partial isometry $v$ in $B(K)$ we have $f(-v) = -f(v)$ and $T_{v} = T_{-v}$. Furthermore, $T_v$ is weak$^*$-continuous and  $f(e) = T_v(e)$ for every minimal partial isometry $e\in (1 - vv^*) B(K)  (1 - v^* v)$;
\item For each minimal partial isometry $v$ in $B(K)$ the equality $f(w) = T_v (w)$ holds for every partial isometry  $w\in (1 - vv^*) B(K)  (1 - v^*v)\backslash\{0\}$;
\item Let $w_1,\ldots,w_n$ be mutually orthogonal non-zero partial isometries in $B(K)$, and let $\lambda_1,\ldots,\lambda_n$ be positive real numbers with $\lambda_1=1$. Then $$f\left(\sum_{j=1}^n \lambda_j w_j\right) = \sum_{j=1}^n \lambda_j f\left(w_j\right);$$
\item For each minimal partial isometry $v$ in $B(K)$ we have $f(x) = T_v (x)$ for every $x\in S(\mathcal{B}_{(1 - vv^*) B(K)  (1 - v^* v)})$;
\item For each partial isometry $w$ in $B(K)$ the element $f(w)$ is a partial isometry;
\item Suppose $v_1,v_2$ are mutually orthogonal minimal partial isometries in $B(K)$ then $T_{v_1} (x) = T_{v_2} (x)$ for every $x\in (1 - v_1v_1^*) B(K)  (1 - v_1v_1^*) \cap (1 - v_2v_2^*) B(K)  (1 - v_2v_2^*)$;
\item Suppose $v_1,v_2$ are mutually orthogonal minimal partial isometries in $B(K)$ then exactly one of the following statements holds:\begin{enumerate}[$(1)$]\item The mappings $T_{v_1}$ and $T_{v_2}$ are complex linear;
\item The mappings $T_{v_1}$ and $T_{v_2}$ are conjugate linear.
\end{enumerate}
\end{enumerate}
\end{theorem}

\begin{proof} $(a)$ Let $v$ be a minimal partial isometry in $B(K)$. Suppose that $$T_v : (1 - vv^*) B(K)  (1 - v^* v)\to  {(1 - f(v)f(v)^*) B(H) (1 - f(v)^*f(v))}$$ is the surjective real linear isometry given by Theorem \ref{t surjective isometries map minimal partial isometries into minimal partial isometries}. Having in mind that $(1 - vv^*) B(K)  (1 - v^*v)\cong B((1 - v^*v)(K), (1 - vv^*)(K))$ and $(1 - f(v)f(v)^*) B(H) (1 - f(v)^*f(v))\cong  B((1 - f(v)^*f(v))(H), (1 - f(v)f(v)^*)(H)) $ are Cartan factors of type 1 and rank $\geq 2$, Proposition 2.6 in \cite{Da} assures that $T_v$ is complex linear or conjugate linear. \smallskip

$(b)$ We keep the notation in $(a)$. Let $T_v, T_{-v} : (1 - vv^*) B(K)  (1 - v^* v)\to  {(1 - f(v)f(v)^*) B(H) (1 - f(v)^*f(v))}$ be the surjective real linear isometries given by Theorem \ref{t surjective isometries map minimal partial isometries into minimal partial isometries}. Lemma 2.5 in \cite{Da} proves that $T_v$ and $T_{-v}$ both are weak$^*$-continuous, while \cite[Theorem 5.1]{ChuDaRuVen} implies that $T_v$ and $T_{-v}$ preserve products of the form $(a,b,c)\mapsto ab^* c + c b^* a$ ($a,b,c\in (1 - vv^*) B(K)  (1 - v^* v)$).\smallskip

Theorem \ref{t surjective isometries map minimal partial isometries into minimal partial isometries} $f(v)$ and $f(-v)$ are minimal partial isometries. By assumptions $\|f(v)-f(-v)\| = \|v +v\| =2$, and hence  by \cite[Lemmas 3.4 and 3.5]{PeTan16} or \cite[Proposition 2.2]{FerPe17} we have $f(-v) = -f(v)$.\smallskip

By the hypothesis on $H$, we can find another minimal partial isometry $e\in (1 - vv^*) B(K)  (1 - v^* v)$. Since $e+v,e-v\in F_e$, applying Theorem \ref{t surjective isometries map minimal partial isometries into minimal partial isometries} we deduce that $$f(v) + T_v(e) = f(e+v) = f(e) + T_e (v)$$ and $$-f(v) + T_{-v} (e) = f(-v) + T_{-v} (e) = f(e-v) = f(e) - T_e (v),$$ where $T_v(e),$ $T_{-v}(e)$ and $T_e (v)$ are minimal partial isometries with $f(v) \perp T_v(e)$, $-f(v)=f(-v) \perp T_{-v} (e)$, and $f(e) \perp T_e (v)$. It follows from Lemma \ref{l minimal partial isometries} above that  $T_e(v) = f(v),$ $T_v (e) = f(e)=T_{-v} (e)$,  and $f(-v) = -T_{e} (v) =-f(v).$\smallskip

We have also shown that $T_v (e) =T_{-v} (e)$ for every minimal partial isometry $e\in (1 - vv^*) B(K)  (1 - v^* v)$. That is, $T_v$ and $T_{-v}$ are surjective complex linear or conjugate linear surjective isometries between Cartan factors of type 1 and rank $\geq 2$. Since $T_v$ and $T_{-v}$ coincide on minimal partial isometries, we deduce by linearity that $T_v$ and $T_{-v}$ both are complex linear or conjugate linear and coincide on finite linear combinations of mutually orthogonal minimal partial isometries.  Finally, since $(1 - vv^*) B(K)  (1 - v^*v)\cong B((1 - v^*v)(K), (1 - vv^*)(K))$ is the weak$^*$-closed span of its minimal
tripotents, we conclude that  $T_v = T_{-v}$.\smallskip

$(c)$ Let $w$ be a non-zero partial isometry in $(1 - vv^*) B(K)  (1 - v^*v)$. Take a minimal partial isometry $w_0$ such that $w= w_0 +(1 - w_0 w_0^*) w  (1 - w_0^* w_0)$. We set $w_0^{\perp} =(1 - w_0 w_0^*) w  (1 - w_0^* w_0)$. Applying Theorem \ref{t surjective isometries map minimal partial isometries into minimal partial isometries} and $(b)$ we get $$f(v) + T_v(w) = f(v + w) = f(v + w_0 +w_0^{\perp}) = f(w_0) + T_{w_0}(v) + T_{w_0} (w_0^{\perp})$$ $$= f(w_0) + f(v) + T_{w_0} (w_0^{\perp}) = f(v) + f(w_0 + w_0^{\perp}) = f(v) + f(w),$$ which proves $(c)$.\smallskip

$(d)$ Let $w_1,\ldots,w_n$ be mutually orthogonal non-zero partial isometries in $B(K)$, and let $\lambda_1,\ldots,\lambda_n$ positive real numbers with $\lambda_1=1$. Pick again a minimal partial isometry $w_0$ such that $w_1= w_0 +(1 - w_0 w_0^*) w_1  (1 - w_0^* w_0)$. Theorem \ref{t surjective isometries map minimal partial isometries into minimal partial isometries} proves that $$ f\left(\sum_{j=1}^n \lambda_j w_j\right) = f\left(w_0 + (1 - w_0 w_0^*) w_1  (1 - w_0^* w_0)+ \sum_{j=2}^n \lambda_j w_j\right) $$ $$= f(w_0) + T_{w_0}\Big( (1 - w_0 w_0^*) w_1  (1 - w_0^* w_0)+ \sum_{j=2}^n \lambda_j w_j \Big) $$ $$= f(w_0) + T_{w_0} ((1 - w_0 w_0^*) w_1  (1 - w_0^* w_0)) + \sum_{j=2}^n \lambda_j T_{w_0}( w_j )  = \hbox{(by $(c)$)}= \sum_{j=1}^n \lambda_j f\left(w_j\right)$$

$(e)$ Since elements in $(1 - vv^*) B(K)  (1 - v^*v)$ can be approximated in norm by finite real linear combinations of mutually orthogonal partial isometries in $(1 - vv^*) B(K)  (1 - v^*v)$ (see \cite[Lemma 3.11]{Horn87}), and $f$  and $T_v$ are isometries, we derive from $(c)$ and $(d)$ that $f(x) = T_v (x)$ for every $x\in S(\mathcal{B}_{(1 - vv^*) B(K)  (1 - v^* v)})$.\smallskip

$(f)$ Let $w$ be a partial isometry in $B(K)$. As before, let $w_0$ be a minimal partial isometry in $B(K)$ satisfying $w= w_0 + w_0^\perp$, with $w_0^{\perp} =(1 - w_0 w_0^*) w  (1 - w_0^* w_0)$. Having in mind that $T_{w_0}$ is a surjective real linear isometry between spaces isometrically isomorphic to $B((1-w_0 w_0^*)(K))$ and $B((1-f(w_0) f(w_0)^*)(H))$, Theorem 5.1 in \cite{ChuDaRuVen} assures that $T_{w_0}$ preserves triple products of the form $\{a,b,c\} =\frac12 (ab^*c + c b^* a),$ and thus $T_{w_0} ( w_0^{\perp})$ is a partial isometry. By Theorem \ref{t surjective isometries map minimal partial isometries into minimal partial isometries}, we get $$f(w) = f(w_0 + w_0^{\perp}) = f(w_0) + T_{w_0} ( w_0^{\perp})$$ is the sum of two orthogonal partial isometries in $B(H)$, and then $f(w)$ is a partial isometry.

$(g)$ Suppose $v_1,v_2$ are mutually orthogonal partial isometries in $B(K)$.  Let us pick a nonzero $x\in (1 - v_1v_1^*) B(K)  (1 - v_1v_1^*) \cap (1 - v_2v_2^*) B(K)  (1 - v_2v_2^*)$. The equality $T_{v_1} (\frac{x}{\|x\|})= f (\frac{x}{\|x\|})= T_{v_2} (\frac{x}{\|x\|})$ holds by $(e)$, and by linearity $T_{v_1} (x) = T_{v_2} (x)$.\smallskip

Finally statement $(h)$ follows straightforwardly from $(a)$ and $(g)$ because the dimensions of $H$ and $K$ are greater than or equal to 3.
\end{proof}

\section{Synthesis of a surjective real linear isometry}

In order to produce a real linear extension of our surjective isometry between $B(H)$ spaces, the next identity principle, which generalizes \cite[Proposition 3.9]{PeTan16}, will play a central role.

\begin{proposition}\label{p surjective isometries extension on projections} Let $H$ and $K$ be complex Hilbert spaces. Suppose that $f: S(B(K)) \to S(B(H))$ is a surjective isometry, and $T:B(K)\to B(H)$ is a weak$^*$-continuous real linear operator such that $f(v)=T(v)$, for every minimal partial isometry $v$ in $B(K)$. Then $T$ and $f$ coincide on $S(B(K))$.
\end{proposition}

\begin{proof} Take a minimal partial isometry $e$ in $B(K)$. By Theorem \ref{t surjective isometries between spheres of B(H)}$(e)$ and the hypothesis $T_e (v) = f(v) = T(v)$ for every minimal partial isometry $v$ in $v\in (1 - ee^*) B(K)  (1 - e^* e)$. Finite real linear combinations of mutually orthogonal minimal partial isometries in $(1 - ee^*) B(K)  (1 - e^* e)$ are weak$^*$-dense in $B(K)$, we therefore deduce from the weak$^*$-continuity of $T_v$ and $T$ that $T_v = T$ on $(1 - ee^*) B(K)  (1 - e^* e)$.\smallskip

Pick a non-zero partial isometry $w$ in $B(K)$, and a minimal partial isometry $w_0$ such that $w=w_0 +(1 - w_0 w_0^*) w  (1 - w_0^* w_0).$ By Theorem \ref{t surjective isometries map minimal partial isometries into minimal partial isometries}, the hypothesis and what we have proved in the first paragraph we obtain $$f(w) = f(w_0) + T_{w_0} ((1 - w_0 w_0^*) w  (1 - w_0^* w_0)) $$ $$= T(w_0)  + T((1 - w_0 w_0^*) w  (1 - w_0^* w_0))  = T(w).$$ We have thus established that $T(w) = f(w)$ for every partial isometry $w$ in $B(K)$. Repeating the arguments in the proof of Theorem \ref{t surjective isometries between spheres of B(H)}$(e)$ we conclude that $T$ and $f$ coincide on $S(B(K))$.
\end{proof}

We have developed enough tools to prove our main result.

\begin{theorem}\label{t Tingley BH spaces} Let $H$ and $K$ be complex Hilbert spaces. Suppose that $f: S(B(K)) \to S(B(H))$ is a surjective isometry. Then there exists a surjective complex linear or conjugate linear surjective isometry $T: B(K) \to B(H)$ satisfying $f(x) = T(x)$, for every $x\in S(B(K))$.
\end{theorem}

\begin{proof} By Riesz's lemma $H$ is finite dimensional if and only if $K$ is. When $H$ and $K$ are finite dimensional, the desired conclusion follows from \cite{Tan2016} or \cite{Tan2016-2}.\smallskip

We assume now that $H$ and $K$ are infinite dimensional. We shall apply the technique in \cite[Theorem 3.13]{PeTan16} to define our real linear isometry. Let $p_1,p_2$ and $p_3$ be three minimal projections in $B(K)$. Given $j\in \{1,2,3\}$, let $T_{p_j} : (1-p_j)B(K)(1-p_j) \to (1-f(p_j)f(p_j)^*)B(H)(1-f(p_j)^*f(p_j))$ denote the surjective real linear isometry given by Theorem \ref{t surjective isometries map minimal partial isometries into minimal partial isometries}.\smallskip

By Theorem \ref{t surjective isometries between spheres of B(H)}$(b)$ and $(h)$, the operators $T_{p_1}$, $T_{p_2}$ and $T_{p_3}$ are weak$^*$-continuous, and they are all complex linear or conjugate linear. We assume that we are in the first case (the second case produces a conjugate linear map).\smallskip

We can mimic the construction done in \cite[Theorem 3.13]{PeTan16} with the appropriate adaptations via the stronger properties developed in Section \ref{Sec:2}. Clearly $B(K)$ admits the following decomposition $$B(K) = \mathbb{C} p_1 \oplus p_1 B(K) p_2\oplus p_2 B(K) p_1 \oplus p_1 B(K) (1-p_1-p_2)$$ $$\oplus (1-p_1-p_2) B(K) p_1
\oplus (1-p_1) B(K) (1-p_1).$$ We define a mapping $T : B(K) \to B(H)$ given by $$T(x) = T_{p_3} (p_1 x p_1) + T_{p_3} (p_1 x p_2 + p_2 x p_1) + T_{p_2}(p_1 x (1-p_1-p_2) +(1-p_1-p_2) x p_1)$$ $$+ T_{p_1}((1-p_1)x (1-p_1)).$$ The mapping $T$ is well defined, complex linear, and weak$^*$-continuous thanks to the uniqueness of the above decomposition and the linearity and weak$^*$-continuity of the mappings $T_{p_1}$, $T_{p_2}$ and $T_{p_3}$ (compare Theorem \ref{t surjective isometries between spheres of B(H)}$(b)$).\smallskip

We shall conclude the proof by applying Proposition \ref{p surjective isometries extension on projections}, for this purpose we shall show that \begin{equation}\label{eq f and T coincide on rank one tripotents second thm} T(e) =f(e), \hbox{ for every minimal partial isometry $v$ in $B(K)$.}
\end{equation} Let $e$ be a minimal partial isometry in $B(K)$. Since dim$(K)=\infty$ there exists a minimal projection $p_4$ satisfying $p_4 \perp p_1, p_2, p_3, e$. The relations of orthogonality imply that $$p_1 e p_1,  p_1 e p_2 + p_2 e p_1, p_1 e (1-p_1-p_2) + (1-p_1-p_2) e p_1, (1-p_1)e (1-p_1)\perp p_4,$$ equivalently the elements $p_1 e p_1,$  $p_1 e p_2 + p_2 e p_1,$ $p_1 e (1-p_1-p_2) + (1-p_1-p_2) e p_1,$ $(1-p_1)e (1-p_1)$ all belong to $(1-p_4) B(K) (1-p_4).$ By definition, Theorem \ref{t surjective isometries between spheres of B(H)}$(g)$ and the previous observation we get $$T(e) = T_{p_3} (p_1 e p_1) + T_{p_3} (p_1 e p_2 + p_2 e p_1) + T_{p_2}(p_1 e (1-p_1-p_2) +(1-p_1-p_2) e p_1)$$ $$+ T_{p_1}((1-p_1)e (1-p_1))= T_{p_4} (p_1 e p_1) + T_{p_4} (p_1 e p_2 + p_2 e p_1)$$ $$ + T_{p_4}(p_1 e (1-p_1-p_2) +(1-p_1-p_2) e p_1)+ T_{p_4}((1-p_1)e (1-p_1))$$ $$= T_{p_4} (e) =\hbox{(Theorem \ref{t surjective isometries between spheres of B(H)}$(e)$)}= f(e) ,$$ which proves \eqref{eq f and T coincide on rank one tripotents second thm} and finishes the arguments.
\end{proof}

We have now tools to extend Theorem \ref{t Tingley BH spaces} for $\ell_{\infty}$-sums of $B(H)$ spaces. In the proof presented here we revise the arguments in the proof of \cite[Theorem 3.12]{PeTan16} and we insert the appropriate modifications.

\begin{theorem}\label{thm Tyngley ellinfty sums} Let $(H_i)_{i\in I}$ and $(K_j)_{j\in J}$ be two families of complex Hilbert spaces. Suppose $f: S\left(\bigoplus_j^{\ell_{\infty}} B(K_j) \right) \to S\left(\bigoplus_i^{\ell_{\infty}} B(H_i) \right)$ is a surjective isometry. Then there exists a surjective real linear isometry $T: S\left(\bigoplus_j^{\ell_{\infty}} B(K_j) \right) \to S\left(\bigoplus_i^{\ell_{\infty}} B(H_i) \right)$ satisfying $T|_{S(E)} = f$.
\end{theorem}

\begin{proof} To simplify the notation, set $A= \bigoplus_j^{\ell_{\infty}} B(K_j)$ and $B=\bigoplus_i^{\ell_{\infty}} B(H_i)$. If $\sharp J \geq 2$, we can pick two different subindexes $j_1$ and $j_2$ in $J$. Let $p_{1}\in B(K_{j_1})$ and $p_{2}\in B(K_{j_2})$ be minimal partial isometries, and let  $T_{p_i} : (1-p_{j_i})A (1-p_{j_i}) \to (1-f(p_{j_i})f(p_{j_i})^*)B (1-f(p_{j_i})^*f(p_{j_i}))$ be the surjective real linear isometry given by Theorem \ref{t surjective isometries map minimal partial isometries into minimal partial isometries}. Let us observe that we can write $A = A_1 \oplus A_2$, where $A_2=\bigoplus_{j\neq j_1}^{\ell_{\infty}} B(K_j)$ and $A_1=B(K_{j_1})$. The symbol $\pi_i$ will stand for the projection of $A$ onto $A_i$. The mapping $T: A\to B$, $T(x) := T_{p_2} (\pi_1 (x)) + T_{p_1} (\pi_2 (x))$ is well defined, real linear and continuous. The operator $T$ is weak$^*$-continuous because $T_{p_1}$ and $T_{p_2}$ are (see Theorem \ref{t surjective isometries between spheres of B(H)}$(b)$). Every minimal partial isometry $v$ in $A$ lies in $A_1$ or in $A_2$. If $v\in A_1$ (respectively, in $A_2$) we have $T(v) = T_{p_2} (v) = f(v)$ (respectively, $T(v) = T_{p_1} (v) = f(v)$) by Theorem \ref{t surjective isometries between spheres of B(H)}$(b)$. Proposition \ref{p surjective isometries extension on projections} assures that $f(x) = T(x)$ for every $x\in S(A)$.\smallskip

If $\sharp I \geq 2$ we can apply the above arguments to $f^{-1}$. We can therefore assume that $\sharp J = \sharp I= 1$ and then the desired statement follows from Theorem \ref{t Tingley BH spaces}.
\end{proof}

\textbf{Acknowledgements} Authors partially supported by the Spanish Ministry of Economy and Competitiveness and European Regional Development Fund project no. MTM2014-58984-P and Junta de Andaluc\'{\i}a grant FQM375.

\end{document}